%% file: attempt_sub.tex
\documentclass{article}
\usepackage{amsmath,amssymb,amsthm,amsfonts,epsfig,fullpage}

\def\Rset{\mathbb{R}}
\def\Nset{\mathbb{N}}

\def\leng{{\mathrm {\tt length}}}
\def\inte{{\mathrm {\tt int}}}
\def\rd{{\mathrm {\tt r\partial}}}
\def\ri{{\mathrm {\tt ri}}}
\def\ah{{\mathrm {\tt ah}}}

\def\cA{{\cal A}}

\def\cH{{\cal H}}

\def\cL{{\cal L}}

\def\cN{{\cal N}}
\def\cHN{{\cal HN}}

\def\cS{{\cal S}}

\def\cX{{\cal X}}

\def\by{\bar{y}}

\def\tx{\tilde{x}}

\def\tpsi{\tilde{\psi}}

\def\tV{\tilde{V}}

\def\neig{{\mathrm {Neighbors}}}   
\def\conv{{\mathrm {conv}}}   
\def\card{{\mathrm {card}}}

\newcommand{\eqd}{\doteq}

\newtheorem{theo}{Theorem}
\newtheorem*{theo*}{Theorem}
\newtheorem{defi}{Definition}

\newtheorem{lemm}{Lemma}

\newtheorem{assu}{Assumption}
\theoremstyle{definition}

\theoremstyle{remark}
\newtheorem{exam}{Example}

\begin{document}

\title{Stability of leaderless multi-agent systems.\\
Extension of a result by Moreau}
\author{David Angeli\thanks{Dipartimento di Sistemi e Informatica, 
University of Florence, Via di S.\ Marta 3, 50139 Firenze, Italy. 
Email: 
\tt angeli@dsi.unifi.it} \and Pierre-Alexandre Bliman\thanks{INRIA, 
Rocquencourt BP105, 78153 Le Chesnay cedex, France. Email: \tt 
pierre-alexandre.bliman@inria.fr}}
\maketitle
\textbf{Abstract.} The paper presents a result which relates 
connectedness of the interaction graphs
in a multi-agent systems with the capability for global convergence to 
a common equilibrium of the system.
In particular we extend a previously known result by Moreau by 
including the possibility of arbitrary bounded
time-delays in the communication channels and relaxing the convexity of 
the allowed regions for the state
transition map of each agent.  
\section{Introduction}
\label{se1}
Recent years have witnessed a growing interest in the study of the 
dynamical
behaviour of the so called multi-agent systems. Roughly speaking these 
can 
be thought of as complex dynamical systems composed by a high number of
simpler units, the agents. Each of them updates its state
according to some rule, whose Input-Output dynamics are typically much 
simpler and much better understood,
and on the basis of the available information coming from the other 
agents.
 All of them, though not necessarily identical, share in fact some 
common feature of interest
(say for instance a given output variable)  and are coupled together by 
communication channels.
The focus of the current research is precisely on how the global 
behaviour of the system,
(for instance questions concerning the global stability or the overall 
synchronization)
is influenced by the topology of the coupling on one hand (this is an 
analysis problem in 
many respects) or the dual question of how to induce a certain desired 
property of the ensemble based on some form 
of local coupling for the agents.
Problems of this nature arise in many different fields, such as in the 
theory of coupled oscillators \cite{JadMot2004,Sepulchre},
in neural networks \cite{Hirsch},
in economics or in the manouvering of groups of vehicles \cite{Naomi}.
For instance in \cite{Lin} the so called \emph{rendezvous} problem is
considered,
namely how to design a local updating rule, based on nearest neighbor
interactions,  which would ensure convergence of all of the agents to 
an unspecified
 common meeting point. Emergence of a global behaviour is therefore a
 consequence of the local updating rule,
 without the need for a leader nor of centralized directions being 
broadcasted.
 
Despite the common traits, the most powerful results are obtained when 
specializing to systems
of a given simple form. Hereby we take a slightly different approach. 
The
emphasis is on how the topology of interconnections between agents 
(possibly
time-varying) affects the convergence of all agents to a common 
equilibrium.
This analysis will be carried out in the presence of 
limited transmission speed of the information between the agents.
In particular, we propose an extension of the contributions by Moreau 
\cite{MOR03,MOR04}, mainly in two directions:
\begin{itemize}
\item
The new setting allows the presence of arbitrary bounded communication 
delays.
\item
A central assumption in the results \cite{MOR03,MOR04}, namely that the 
future evolution of the studied system is constrained to occur in the 
convex hull of the agents states, is removed.
\end{itemize}
The first aspect comes as a very natural question both from a practical 
and a theoretical
point of view. Communication delays are in fact ubiquitous in the 
``real'' world and it is
well-known their potential destabilizing effect in conjunction with 
feedback
loops, here induced by the graph topology of the communication channels 
which need not be of a hierarchical type.
 It is therefore remarkable to see how, at least in the specific set-up 
we are considering,
this destabilizing effect does not take place and the same global 
behaviour of the multi-agent
system in terms of convergence to a common equilibrium follows also in 
the extended set-up.                  \\
\newline
The second extension deals with convexity issues; one of the technical 
tools used in order to enforce a common
behaviour in systems whose state takes value in Euclidean space, is to 
have local evolutions point always
inside the convex hull of all variables. This makes life easier in a 
certain respect but it is an unnatural assumption
in more general contexts, for instance when oscillators networks are 
considered (these are typically modeled as systems evolving on a torus)
or systems evolving in partially obstructed Euclidean spaces (for 
instance on a plane minus a circle).
 Relaxing convexity is meant as a first step in the quest for stability 
conditions which can work in more general spaces.

Before going on further, we present the main elements of this 
construction, developed below.
The multi-agent system under study will be described by a {\em 
time-dependent graph} $\cA(t)$, describing the transfer of information 
between 
the agents at time $t$, and a {\em set of rules\/} according to which 
each agent updates its state at time $t+1$.
The definition of the latter is done by the introduction of two types 
of objects, which we present now 
(complete definitions are to be found in Section \ref{se2} below).

\begin{itemize}
\item
A {\em set-valued map\/} $\sigma$, is defined, which associates to the 
set of present and past states of the agents a compact set in the state 
space common to all the agents.
This map will play the central role of a {\em set-valued Lyapunov 
function\/} for the system.

\item
It is then necessary to define the rules according to which the agents 
update their state, given the (possibly delayed) information on the 
position of the other agents they received.
For this, each agent $k$ is attributed a set-valued map $e_k$ which, 
{\em given the communication graph $\cA(t)$}, defines the set of 
allowed 
positions $e_k(\cA(t))$.
An important point here is that, whatever the information received by 
each agent, the new positions cannot induce an increase of the 
set-valued Lyapunov function along the trajectories.

\end{itemize}
The definition of the new class of multi-agent systems studied here is 
done and commented in Section \ref{se2}.
The stability is studied and the main results are given in Section 
\ref{se3}.

\paragraph{Notations}
As often as possible, we use the notations introduced by Moreau 
\cite{MOR03,MOR04}.
Following him, we distinguish between the inclusion, denoted 
$\subseteq$, and the strict inclusion, denoted $\subset$.
The topological interior of a set is denoted $\inte$.

We study systems with $n$ agents whose position at time $t$ are written 
as $x_1(t),\dots, x_n(t)$ in the finite-dimensional space $X$.
In the setting introduced in Moreau's contributions, the corresponding 
overall state variable is $x(t)=(x_1(t),\dots, x_n(t)) \in X^n$.
Here, we consider systems with delay smaller than a given integer 
$h>0$.
In consequence, the complete state variable of the system is 
$(x_1(t),x_1(t-1),\dots, x_1(t-h+1),\dots, x_n(t), \dots, 
x_n(t-h+1))\in X^{hn}$.

We denote $\tx =(x_1,\dots, x_{hn})$ an arbitrary element of $X^{hn}$ 
and, when considering the dynamical system, we write 
$\tx_k(t)=(x_k(t),x_k(t-1),\dots, x_k(t-h+1))$ for all 
$k\in\cN\eqd\{1,\dots, n\}$ and 
$\tx(t)=(\tx_1(t),\dots, \tx_n(t))$.
We also use the corresponding decomposition of any element $\tx$ of 
$X^{hn}$ as
$\tx=(\tx_1,\dots, \tx_n)$ (which amounts to identify $X^{hn}$ to 
$(X^h)^n$).
When needed, any $\tx_k\in X^h$ is decomposed according to 
$\tx_k=(x_{k,0},\dots, x_{k,h-1})$,
in such a way that for the variables of the dynamical systems under 
study $x_{k,j}(t)=x_k(t-j)$, $k\in\cN$, $j\in\cH\eqd\{0,\dots, h-1\}$. 
Similarly we denote by $\cHN \eqd \{ 1,2, \ldots, hn \}$.
The previous notation is necessary, in order to be able to distinguish 
between the delayed and the actual values of the position of the 
agents.
Coherently with the notations introduced above, we sometimes abbreviate 
$x_{k,0}$ and write simply $x_k$.

Last, given any $\tilde{x} \in X^{hn}$ we often need to embed it on $2^X$,
according to the following rule: $\pi ( \tilde{x} ) \eqd$ $\{ x_1, x_2,
\ldots x_{nh} \}$. In this way the state of the system is mapped to a
finite collection of points in the $X$ space.\\

Finally, for the comfort of the reader we indicate that the Theorems 1, 2, 
3 and 5 in reference \cite{MOR03} are numbered respectively  4, 1, 2 
and 5 in \cite{MOR04a}.

\section{A class of multi-agent dynamical systems}
\label{se2}

This section is devoted to the presentation of the dynamical system 
under study.
We study here a special class of {\em nonlinear difference inclusions 
with delay},
that we write:
\begin{equation}
\label{equa}
x_k(t+1) \in  e_k(\cA(t))(\tx (t) )\ .
\end{equation}

Recall that $x_k(t)$ represents the ``position" at time $t$ of the 
agent $k$.
The evolution of the latter depends upon the complete system state 
$\tx(t)$ (including delayed components), through the time-varying  map 
$e_k(\cA(t))$.
For a trajectory of \eqref{equa}, we call {\em decision set of agent 
$k$} at time $t$ the value taken by $e_k(\cA(t))(\tx (t))$.
The specificity of the problem lies in these maps: they depend upon 
the topology of the inter-agent communications, modeled by the {\em 
graph\/} $\cA(t)$.

The modeling of the communication network is presented below in Section 
\ref{se21}.
The construction of the decision sets inside which, given the 
communication network, each agent may update the value of its state, is 
made 
in Section \ref{se22}.
Last, we provide some examples in Section \ref{se23}.

\subsection{Inter-agent communications modeling}
\label{se21}

The first ingredient of the construction is the family of {\em 
continuous 
set-valued maps} $e_k(\cA) : X^{hn}\rightrightarrows X$ taking on 
{\em compact 
values}, and defined for $k\in\cN$ and any directed graph $\cA$.
 The latter will define, according to the position of the other agents,
 in which subset of $X$ agent $k$ 
is allowed to choose its future state.

 Here, we are concerned by information transfer from the past to the 
present.
 In other words, we need to consider graphs in $X^{hn}$ 
linking some past and/or present values $x_k(t-j)$ of the states of an 
agent $k$ to 
another agent $l$.
This motivates the following adaptations of some notions of the graph 
theory.
Recall that $\cN=\{1,\dots, n\}$, $\cH=\{0,\dots, h-1\}$, where $n$ is 
the number of agents and $h-1$ the larger transmission delay.
 
 \begin{defi}
 \label{de1}
 We call {\em directed graph with delays} any subset $\cA$ of 
$(\cN\times\cH)\times \cN$ satisfying $((k,0),k)\not\in\cA$ for all 
$k\in\cN$.
 The elements of $\cN$ are called {\em nodes} and an element 
$((k,j),l)$ of $\cA$ is referred to as an {\em arc from $k$ to  $l$}.
 A node  $k\in\cN$ is said {\em connected} to a node  $l\in\cN$ if 
there is a path from $k$ to $l$ in the graph with delays which respects 
the 
orientation of the arcs.
Given a sequence of directed graphs with delays $\cA(t)$, a node  
$k\in\cN$ is said {\em connected\/}
to a node  $l\in\cN$ {\em across an interval $I\subseteq\Nset$} if 
 $k$ is connected to $l$ for the graph with delays $\bigcup_{t\in 
I}\cA(t)$.
 \end{defi}
 
\begin{figure}
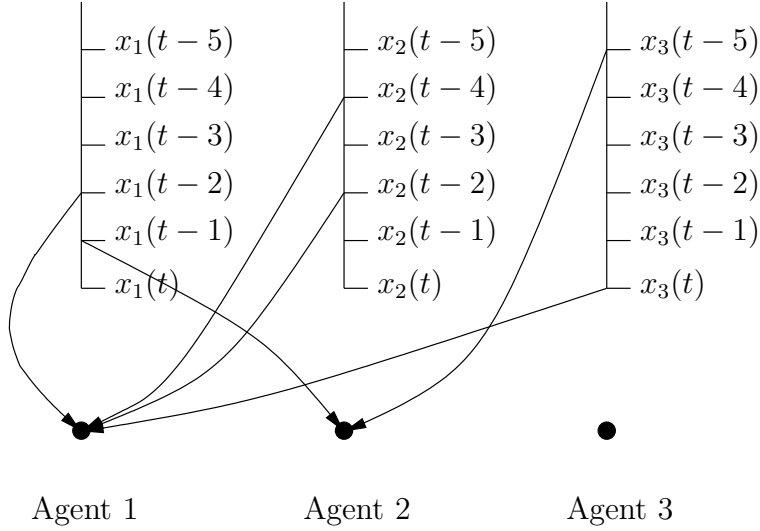

 \begin{center}
\input fig7.pstex_t
  \end{center} 
 \caption{An example of graph with delays for a system with three 
agents.}
   \label{fi7}
 \end{figure}

Figure \ref{fi7} provides an example of graph with delays.
For the graph represented therein, agents 1 and 2 are mutually connected and agent 3 is connected to 1 and 2, but neither 1 nor 2 is connected to 3.
 Notice that generally speaking there may exist more than one arc between two distinct 
nodes, and that a node may be connected to itself
(via delayed values).

 \begin{defi}
 \label{de2}
Consider a directed graph with delays $\cA$ and a nonempty subset 
$\cL\subseteq\cN$.
The set $\neig(\cL,\cA)$ is the set of those nodes 
$k\in\cN\setminus\cL$ for which there is $l\in\cL$ such that (at least) 
one arc from $k$ to 
$l$ exists.
When $\cL$ is a singleton $\{l\}$, the notation $\neig(l,\cA)$ is used 
instead of $\neig(\{l\},\cA)$.
 \end{defi}
 
 We impose to the maps $e_k$ the following assumption.
 
 \renewcommand{\theassu}{\Alph{assu}}

 \begin{assu}
 \label{as1}
 For all $k\in\cN$ and all directed graph with delays $\cA$,
the set-valued map $e_k$ is continuous and takes on compact values.
Moreover,
\begin{itemize}
  \item
  $e_k(\cA)(\tx) =\{x_k\}$ if\;
$\{x_{i,j}\ :\ ((i,j),k)\in\cA\} 
=\{x_k\}$; 
  
  \item
 $e_k(\cA)(\tx) \subset\ri\ \sigma\left(
\{x_k\} \cup\{x_{i,j}\ :\ ((i,j),k)\in\cA\}
 \right)$ otherwise.
\end{itemize}

\end{assu}
The exact meaning and the properties of the set-valued $\ri\ \sigma$ 
are the subject of Section \ref{se22}.
However, we may already make some remarks on the form of the right-hand
side of the problem.
Clearly, Assumption \ref{as1} implies that the evolution of each agent 
depends only upon the possibly delayed information received from its 
neighbors.
The case where
$\{x_{i,j}\ :\ 
((i,j),k)\in\cA\}=\{x_k\}$ is realized when either the agent $k$ has no 
neighbor and the set involved in the formula is empty, or 
all the (possibly delayed) positions received from the neighboring 
agents are also 
equal to the present position $x_k$ of agent $k$;
in this case, no motion is allowed.
We shall see below that in the present framework the use by each agent 
of the present value of its own position is mandatory for stability, 
see counterexample in Example \ref{ex4}.
 
\subsection{Construction of the decision sets}
\label{se22}

The second ingredient necessary for the construction of the dynamical 
system under study is a {\em  set-valued map} $\sigma : 
2^X \rightrightarrows X$, taking on {\em compact values}.
It has a central role in the definition of the dynamics, and it will be 
shown afterwards (cf.\ in particular the proof of Theorem 
\ref{th2}) that it plays the role of a ``set-valued Lyapunov function'' 
for the studied system.

In order to state the properties that $\sigma$ should fulfil, we have 
to introduce beforehand some notions.
First of all, define $\cS$, a set of subsets of $X$ in which $\sigma$ 
will be compelled to take on its values, as:
\begin{equation}
\label{cS}
\cS\eqd\{
S\subset X\ :\ S \text{ compact and $\exists \varphi : X\to X, \varphi$ 
bijective, $\varphi, \varphi^{-1}$ Lipschitz and 
$\varphi(S)$ convex}
\} \ .
\end{equation}

Important 
consequences will proceed from 
the fact that $\sigma$ takes on values in $\cS$, inherited from 
properties of $\cS$ summarized in the 
following 
result (see proof of Lemma \ref{le15} in Appendix).

\begin{lemm}
\label{le15}
Let $\cS$ be defined by \eqref{cS}.
\begin{enumerate}
\item
for any $S\in\cS$, the function $d_S (x^0,x^1):S\times S\to 
[0,+\infty)$ defined as 
\[
d_S(x^0,x^1)\mapsto \inf\left\{
\leng(\psi)\ :\ \psi: [0,1] 
\stackrel{\textrm{Lipschitz}}{\longrightarrow}S, 
\psi(0)=x^0, \psi(1)=x^1
\right\}
\]
is well-defined and continuous.
Define $\mu : \cS\to \Rset^+$ by:
\begin{equation}
\label{mu}
\mu(S)\eqd\max_{x^0, x^1\in S} d_S (x^0,x^1).
\end{equation}
Then, for all $S\in\cS$,
\begin{itemize}
\item
$\mu(S)<+\infty$.
\item
$\mu(S)=0$ if and only if $S$ is a singleton.
\item
$\mu(S)$ is at least equal to the (euclidian) diameter of $S$, and 
equal to this value if $S$ is convex.
\item
$\mu$ is lower semicontinuous in $S$,
but nowhere continuous.
\end{itemize}

\item
for any $S\in\cS$, let $\varphi$ be as in \eqref{cS} and
\[
\ri (S) \eqd \varphi^{-1}\left(
\ri(\varphi(S))
\right)\ ,
\]
where $\ri(\varphi(S))$ designates the relative interior of the convex 
set $\varphi(S)$, i.e.\ its interior when regarded as a topological 
subspace of its affine hull.
Then, for all $S\in\cS$,
\begin{itemize}
\item
$\ri(S)$ is independent of the choice of $\varphi$.
\item
$\ri(S)=\emptyset$ if and only if $S$ is a singleton.
\item
$\inte\ S\subseteq\ri\ S\subset S$.
\item
$\ri(S)$ is the relative interior of $S$ if $S$ is convex.

\end{itemize}
\end{enumerate}
\end{lemm}

Lemma \ref{le15} permits to measure the distance between points of a 
set $S\in\cS$ ``along the arcs".
It permits to define extended notions of diameter and of relative 
interior, which coincide with the usual ones for convex subsets of $X$.
By definition, we call ``relative boundary" of sets $S$ in $\cS$ the 
following set:
\[
\rd (S) \eqd S\setminus\ri (S)\ .
\]
Also, according to the definition of $d_S$ in Lemma \ref{le15}, we define, for any subsets $S', S''$ of a set $S$ in $\cS$ the {\em $S$-distance from $S'$ to $S''$} as:
\begin{equation}
\label{last}
d_S(S',S'') \eqd \inf_{x^0\in S', x^1\in S''} d_S(x^0,x^1)\ .
\end{equation}

We now gather the properties that $\sigma$ must fulfil,
 and afterwards 
comment on their meaning and consequences.

\begin{assu}
\label{as2}
The set-valued map $\sigma:2^X \to\cS$ is continuous with respect to the topology 
induced by Hausdorff metric. 
Moreover,
the following should hold:

\begin{enumerate}

\item
\label{assumption0}
$S\subseteq\sigma(S)$ with equality if $S$ is a singleton.

\item
\label{assumption1}
$\sigma(S)= \sigma \circ \sigma(S)$ for all $S \in 2^X$.

\item
\label{assumption2}
$S'\subseteq S \; \Rightarrow \; \sigma(S')\subseteq\sigma(S)$ for all 
$S, S' \in 2^X$.

\item
\label{assumption3}
If $S$ is not a singleton, for all $x\in S$, there exists 
$\Sigma_x\subseteq\rd\sigma(S)$ such that 
$\Sigma_x\cap S\neq\emptyset$ and $x\not\in\Sigma_x$. Moreover, if 
$S'\subseteq\sigma(S)$:
\begin{enumerate}
\item
\label{assumption3a}
if $\ri\ \sigma(S')\cap\Sigma_x\neq\emptyset$, then $S'\subseteq\Sigma_x$ (and in particular, $x\notin S'$).

\item
\label{assumption3b}
if $d_S(S',\Sigma_x)>0$, then
$\mu(\sigma(S'))<\mu(\sigma(S))$.

\end{enumerate}
\item
\label{assumption4}
$\mu\circ\sigma$ is continuous.

\end{enumerate}
\end{assu}

Remark that at this point, the problem under study is fully 
understandable:
our goal is to find stability conditions for systems defined by 
\eqref{equa}, where the maps $e_k$ verify Assumption \ref{as1} for a 
given map 
$\sigma$ fulfilling Assumption \ref{as2}, and where the meaning of the 
relative interior $\ri$ has been defined previously by Lemma 
\ref{le15}.\\

Important consequences of Assumptions \ref{as2}.\ref{assumption0} to 
\ref{as2}.\ref{assumption4} are now 
stated.
We shall see further (see Theorem \ref{th1}) that Assumptions 
\ref{as2}.\ref{assumption0}--\ref{as2}.\ref{assumption2} are indeed 
sufficient to forbid 
increase along time of the natural 
set-valued Lyapunov function of the system.
The additional Assumptions 
\ref{as2}.\ref{assumption3}--\ref{as2}.\ref{assumption4} 
induce the {\em strict\/} decrease of the set-valued Lyapunov function 
(see Theorem \ref{th2}) .

We provide in the following lemma a direct consequence of Assumption 
\ref{as2}.\ref{assumption0}.
\begin{lemm}
\label{le16}
Assume Assumption \ref{as2}.\ref{assumption0} is 
fulfilled.
Then, for any $S \subset X$ we have: $\card\ S>1 \Rightarrow \ri\ 
\sigma(S) \neq \emptyset$ and
$\mu ( \sigma(S))>0$.
\end{lemm}
\begin{proof}
By Assumption \ref{as2}.\ref{assumption0} it follows from $\card\ S>1$ 
that $\card\ \sigma 
(S) >1$. Since $\sigma(S)$ is not a singleton,
it is homeomorphic to a sphere of non-zero dimension.
Hence, $\ri\ \sigma (S) \neq \emptyset$. Moreover, $\mu ( \sigma (S)) 
\geq 
\textrm{diam} ( \sigma (S)) >0$, since the euclidean diameter of a set 
vanishes if and only if the set is a singleton. 
\end{proof}
 
 We now come to the central hypothesis, stated in Assumption 
\ref{as2}.\ref{assumption3}.
This Assumption applies to arbitrary (but non trivial) groups of agents 
$S$, which may comprise indifferently 
true agents or ``virtual'' agents, viz. informations relative
to the position of a true agent at previous sampling times.  
More closely, for each agent $x$, there exists a portion of the 
boundary of $\sigma(S)$, denoted by $\Sigma_x$,
whose elements are {\em irreversibly\/} attracted outside of it when 
using information received
 from any agent not in $\Sigma_x$  (such as $x$ itself) according to 
the rule edicted in Assumption \ref{as1}.
The second part of Assumption \ref{as2}.\ref{assumption3} imposes that 
such an irreversible 
escape from $\Sigma_x$ comes with a {\em strict decrease\/} 
of the diameter of the set-valued Lyapunov function of the system
(for {\em convex\/} sets $S,S'\subseteq X$, $S'\subseteq S$ 
implies 
$\mu(S')\leq \mu(S)$, but this is not true for general sets in $\cS$ 
defined by \eqref{cS}).

Generally speaking, the set $\Sigma_x$ looks like as an union of 
``faces" of $\rd\sigma(\tx)$ containing an extremity of each geodesic
 of maximal length in $\sigma(S)$ (that is of length equal to 
$\mu(\sigma(S))$), except the point $x$ itself if it belongs to 
 $\rd\sigma(\tx)$.
Remark that sets $\Sigma_x, \Sigma_y$ associated to different points
$x,y$ in $S$ may be equal, and that the union of the sets $\Sigma_x$ over 
$x\in S$ 
does not have to cover $\rd \sigma(S)$.\\

Similarly to what happens within Moreau's setting, one has the 
following result.
\begin{lemm}
\label{le17}
Assume Assumptions 
\ref{as2}.\ref{assumption0}-\ref{as2}.\ref{assumption3}  be fulfilled.
Then, $\card ( S ) > 1 \; \Rightarrow \; \card ( S \cap 
\rd \sigma(S) ) \geq 2 $.
\end{lemm}
\begin{proof}
Let $x \in S$. Since $S$ is not a singleton, by Assumption 
\ref{as2}.\ref{assumption3}, there 
exists $\Sigma_x \subset \rd \sigma(S)$ so that $\Sigma_x \cap S \neq
\emptyset$.
Let $y \in \Sigma_x \cap S$. Since $y \in S$ we can apply the same 
property to $y$ in order to conclude that there exists
$\Sigma_y \subset \rd \sigma (S)$ and such that $\Sigma_y \cap S \neq 
\emptyset$. Let $z \in \Sigma_y \cap S$. Since by assumption 
$y \notin \Sigma_y$ we have $y \neq z$. Therefore $y$ and $z$ both 
belong to $S \cap \rd \sigma(S)$ and are different from each other.
This completes the proof of the lemma. 
\end{proof}

Last, notice that Lemma \ref{le15} and the continuity assumption on 
$\sigma$ implies that the map $\mu\circ\tV$ is already lower 
semicontinuous on $X^{hn}$.
Assumption \ref{as2}.\ref{assumption4} thus represents a slightly 
stronger 
regularity assumption.

\subsection{Examples}
\label{se23}

We present here different examples and counter-examples of maps 
$\sigma$ fulfilling the properties previously defined.

 \begin{exam}[convex hull]
 \label{ex1}
 In Moreau's work, $\sigma(S)$ is taken to be the convex hull of $S$,
see Figure \ref{fi1}.
 One may check easily that Assumptions \ref{as2}.\ref{assumption0} to 
\ref{as2}.\ref{assumption4} are all
fulfilled.
Here, the sets $\Sigma_x$ involved in Assumption 
\ref{as2}.\ref{assumption3} can be defined as 
follows:
\[ \Sigma_x \eqd \bigcup_{c \in TC_{\sigma(S)}(x), |c|=1 } x + \max \{ t: 
x+ ct \in \sigma(S) \} c\ ,
\]
where $TC_{\sigma(S)}(x)$ denotes the Bouligand contingent cone to the set $\sigma(S)$ at $x$
(otherwise called tangent cone, as  $\sigma(S)$ is convex here; see \cite[ pp.\ 176--177 and 219]{aubin} for details).
  \begin{figure}
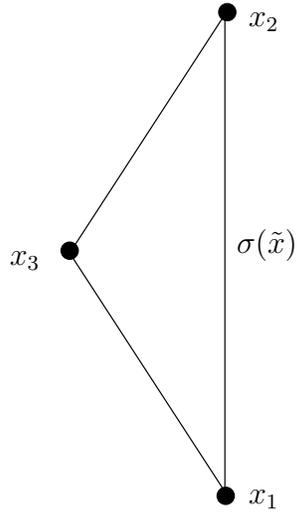

 \begin{center}
\input fig1.pstex_t
  \end{center} 
 \caption{The convex-hull, Moreau's set-valued Lyapunov function.}
  \label{fi1}
 \end{figure}
  \end{exam}
  
 \begin{exam}[a different convex example]
 \label{ex2}
 For a given basis $e_j$, $j=1,\dots,p$ of $X$, take
 \[
 \sigma(S)\eqd  \left[
 \min_{x \in S} e_1^T x,\max_{x \in S} 
e_1^T x
 \right]
 \times\dots\times
  \left[
 \min_{ x \in S} e_p^T x,\max_{x \in S} 
e_p^Tx_i
 \right] \ .
 \]
In this example, the convex hull is applied ``componentwise", see 
Figure \ref{fi2}.
 Remark that $\conv(S)\subseteq\sigma(S)$ for this case,
but this relation is not mandatory, see Example \ref{ex3}  below.

In the example depicted on Figure \ref{fi2}, one may check that the 
choice consisting in taking for 
$\Sigma_x\eqd$ $\bigcup_{c \in TC_{\sigma(S)}(x), |c|=1 } x + \max \{ t: x+ 
ct \in \sigma(S) \} c$, fulfills the Assumptions.

  \begin{figure}
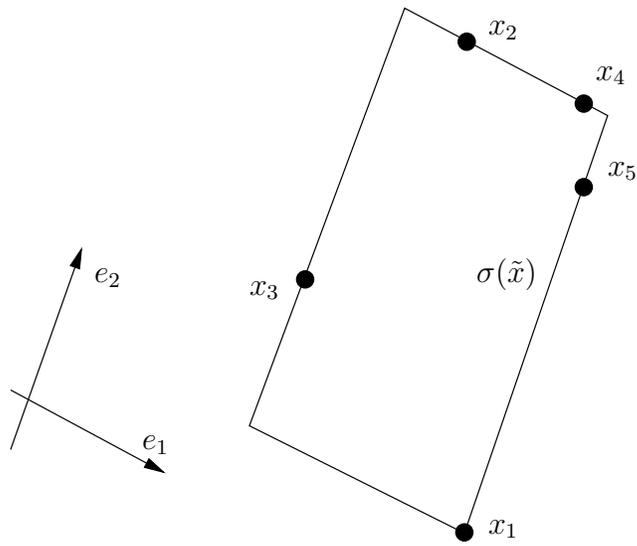

 \begin{center}
\input fig2.pstex_t
  \end{center} 
 \caption{Illustration of Example \ref{ex2}.}
  \label{fi2}
 \end{figure}
  \end{exam}

\begin{figure}
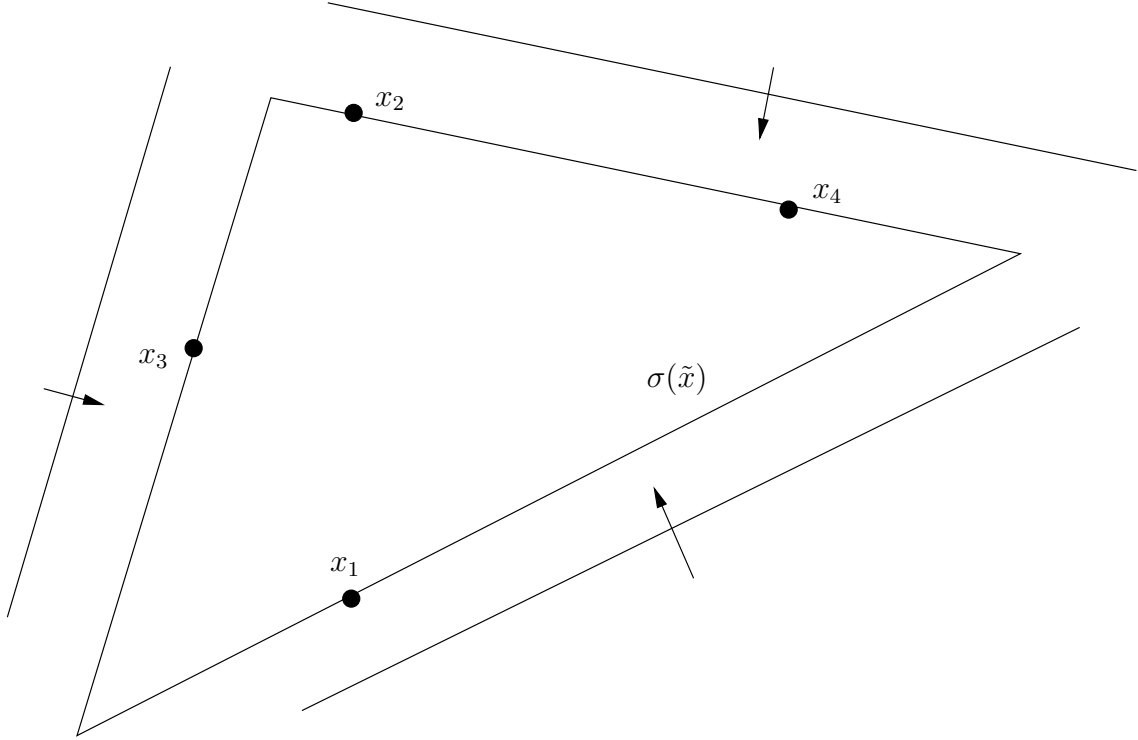

 \begin{center}
\input fig3.pstex_t
  \end{center} 
 \caption{Other convex examples of set-valued Lyapunov function, see 
Example \ref{ex25}.}
    \label{fi3}
 \end{figure}

 \begin{exam}[other convex examples]
 \label{ex25}
 One may also define $\sigma(S)$ as the smaller set containing $S$ and 
with boundary parallel to 
given 
$p+1$ non-parallel hyperplans (where $X=\Rset^p$), see Figure 
\ref{fi3}.
More precisely, let $\Sigma= \conv (S) $ and $e_1, 
\dots, e_{p+1}$ be $(p+1)$ vectors in $X$ 
such that for some positive $\lambda \in \Rset^{p+1}$ we have $\sum_j 
\lambda_j e_j= 0$.
 The set $\sigma(S)$ is a polytope defined as:
 $\left\{
 x\in X\ :\ e_j^Tx \leq \max_{x' \in \Sigma} e_j^{T} x',\ j=1,\dots, 
p+1
 \right\}$,  
 containing the points $x_1,\dots, x_{hn}$.
 Symmetrically we may define
 $\sigma(x)= \left\{
 x\in X\ :\ e_j^Tx \geq \min_{x' \in \Sigma} e_j^{T} x',\ j=1,\dots, 
p+1
 \right\}$. 
Similarly to what occurs in Example \ref{ex2}, one may take for 
$\Sigma_x$ the portion of the boundary obtained by following the 
vectors coming out from the tangent cone at $x$
all the way to their extreme intersection point with the boundary of 
$\sigma(S)$, 
and the Assumptions 
\ref{as2}.\ref{assumption0}-\ref{as2}.\ref{assumption4} are fulfilled.
\end{exam} 
 
 Remark that the smallest ball or the smallest hypercube containing $S$ 
does {\em not\/} fulfil the requested 
properties.
 For instance the smallest circle containing a triangle never contains 
the smallest circle containing the shortest of its edges,
 which violates monotonicity of the map $\sigma$.
 
 \begin{exam}[nonconvex examples]
 \label{ex3}
 For any bijective transformation $\varphi : X\rightarrow X$ which is  
Lipschitz together with its inverse, one may take
 \[
 \sigma_\varphi(S)\eqd
 \varphi^{-1}\left(
 \sigma( \varphi(S) )
 \right)\ ,
 \]
 where $\sigma$ fulfils all the Assumptions.
 In general $\sigma_\varphi(S)\not\subseteq \conv (S)$ 
and is not convex: indeed, this latter property is not essential.
 Such an example of nonconvex sets is given in Figure \ref{fi5}, 
obtained for
 $X=\Rset^2$, $x_1 = \begin{pmatrix} 2 \\ 0 \end{pmatrix}$,
 $x_2=\begin{pmatrix}
 1\\5
 \end{pmatrix}$, 
 $x_3=\begin{pmatrix}
 0\\-1
 \end{pmatrix}$,
 $\varphi(x)=\begin{pmatrix}
 \cos \alpha\|x\|^2&\sin  \alpha\|x\|^2 \\ -\sin  
\alpha\|x\|^2 &\cos  \alpha\|x\|^2
 \end{pmatrix} x$, $\alpha=0.04$, and $\sigma(S)=\conv (S)$.
 
 Notice that, generally speaking,  the systems generated along this 
principle are such that the map $\varphi$ in \eqref{cS} is identical 
for {\em all\/} the sets $\sigma(S)$.
The sets $\Sigma_x$ may be obtained as for Example \ref{ex1},
up to transformation by $\varphi$.
 
  \begin{figure}
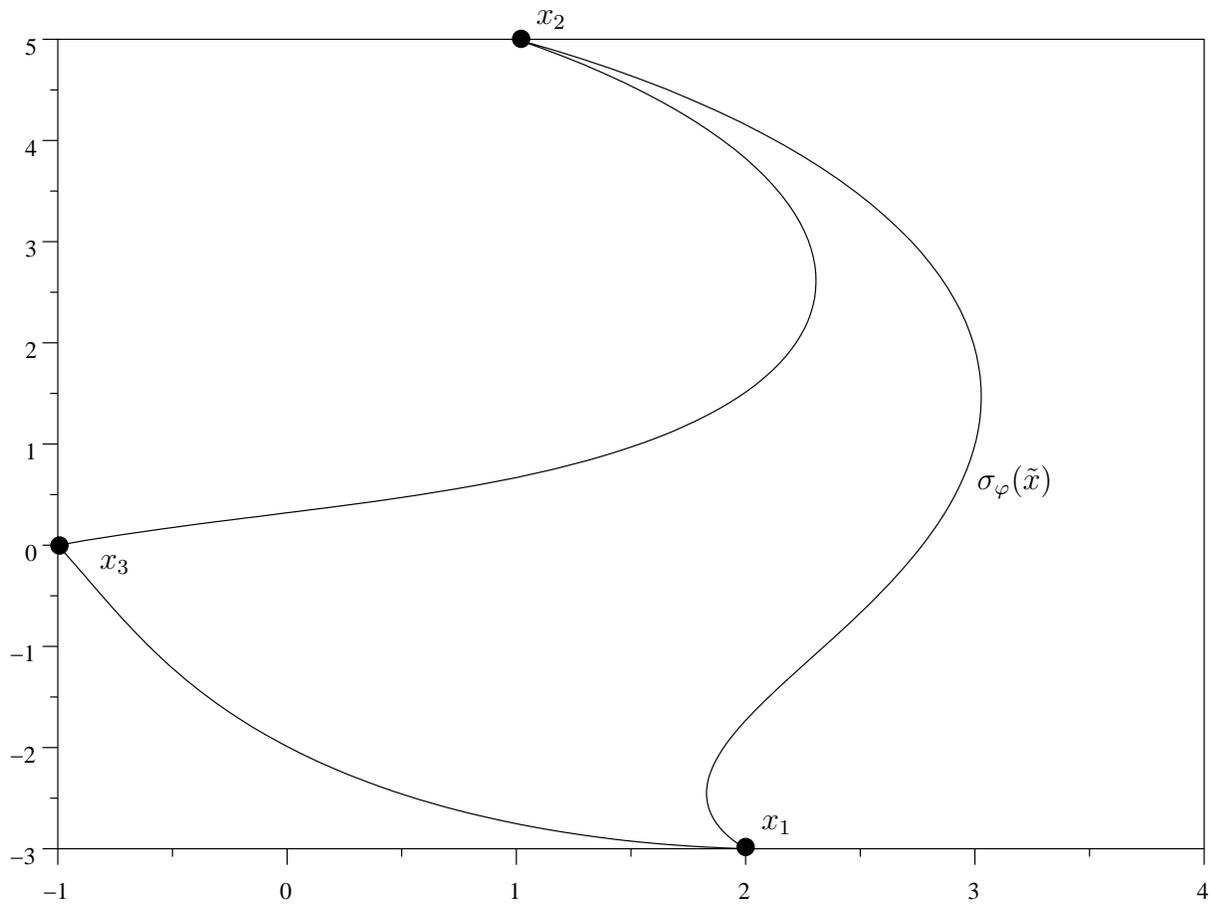

\hspace{-2cm}
\input fig5.pstex_t
 \caption{An example of map $\sigma$ giving rise to nonconvex sets.
 Notice that $\conv (S) \not\subseteq\sigma ( S )$, and that 
$\mu(\sigma(S))$ is larger than the diameter of $\conv (S)$.}
   \label{fi5}
 \end{figure}
 \end{exam}
 
 \begin{exam}[intersection of decision sets]
 \label{ex35}
 When $\sigma, \sigma'$ fulfil the properties stated above, an 
interesting issue is to see whether $\sigma\cap\sigma'$ do.
One verifies easily that Assumptions 
\ref{as2}.\ref{assumption0}--\ref{as2}.\ref{assumption2} are 
automatically 
fulfilled.
The validity of \ref{as2}.\ref{assumption3} and \ref{as2}.\ref{assumption4} depends upon the 
configuration of the sets $\Sigma_x$, $\Sigma_x'$ corresponding to $\sigma$ and 
$\sigma'$.
In Figure \ref{fi6} is presented an example where the resulting map 
fulfills all the properties.\\

  \begin{figure}
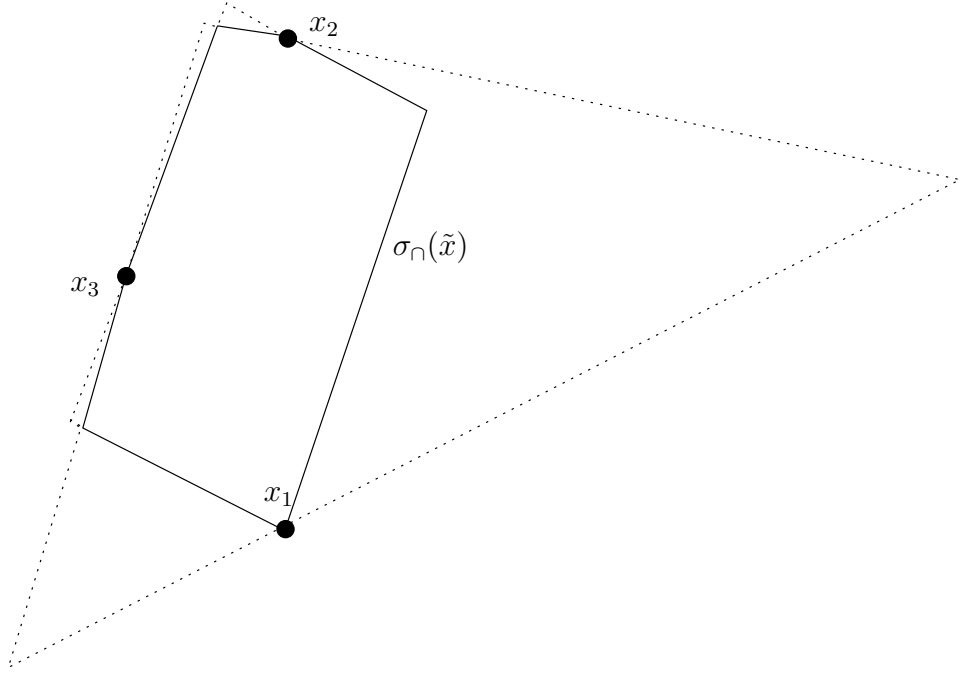

 \begin{center}
\input fig6.pstex_t
  \end{center} 
 \caption{Map obtained by intersection of the 
maps from Figures \ref{fi2} and \ref{fi3}.}
   \label{fi6}
 \end{figure}
 \end{exam}
 
\section{Results}
\label{se3}

Before stating the results of this paper, we recall the notions under 
discussion below, see \cite{MOR03,MOR04a}.
As in Moreau's papers, we call {\em equilibrium point\/} any element of the state space which
is the constant value of an {\em equilibrium solution}.

\begin{defi}
Let $\cX$ be a finite-dimensional Euclidean space and consider a 
continuous set-valued map $e\ :\ \Nset\times\cX\rightrightarrows \cX$ taking on closed values, giving rise to the difference inclusion
\begin{equation}
\label{eque}
x(t+1) \in e(t,x(t))\ .
\end{equation}
Consider a collection of equilibrium solutions of this equation and 
denote the corresponding set of equilibrium points by $\Phi$.
By definition, $\varphi\in\Phi$ if and only if $\varphi\in e(t,\varphi)$ for all $t\in\Nset$.

With respect to the considered collection of equilibrium solutions, the 
dynamical system is called

\begin{enumerate}
\item
{\em stable} if for each $\varphi\in\Phi$, for all $c_2>0$ and for 
all $t_0\in\Nset$, there is $c_1>0$ such that every solution $\zeta$ of \eqref{eque}
satisfies: if $|\zeta(t_0)-\varphi|<c_1$ then $|\zeta(t)-\varphi|<c_2$ for all $t\geq t_0$.

\item
{\em bounded} if for each  $\varphi\in\Phi$, for all $c_1>0$ and for 
all $t_0\in\Nset$, there is $c_2>0$ such that every solution $\zeta$ of \eqref{eque}
satisfies: if $|\zeta(t_0)-\varphi|<c_1$ then $|\zeta(t)-\varphi|<c_2$ for all $t\geq t_0$.

\item
{\em globally attractive} if  for each  $\varphi_1\in\Phi$, for all 
$c_1,c_2>0$ and for all $t_0\in\Nset$, there is $T\geq 0$ such that every 
solution $\zeta$ of \eqref{eque} satisfies: if $|\zeta(t_0)-\varphi_1|<c_1$ then there 
is $\varphi_2\in\Phi$ such that $|\zeta(t)-\varphi_2|<c_2$ for all 
$t\geq t_0+T$.

\item
{\em globally asymptotically stable} if it is stable, bounded and 
globally attractive.

\end{enumerate}
If $c_1$ (respectively $c_2$ and $T$) may be chosen independently of 
$t_0$ in Item 1 (respectively Items 2 and 3) then the dynamical system is 
called uniformly stable (respectively uniformly bounded and uniformly 
globally attractive) with respect to the considered collection of 
equilibrium solutions.
\end{defi}
Notice that the above notions are uniform with respect to all trajectories of \eqref{eque}.

We now state a first result on boundedness and (simple) stability, 
analogous to \cite[Theorem 2]{MOR03}.

\begin{theo}
\label{th1}
Assume that Assumptions \ref{as1} and 
\ref{as2}.\ref{assumption0}--\ref{as2}.\ref{assumption2} are fulfilled.
Then the discrete-time system \eqref{equa} is uniformly globally 
bounded and uniformly globally stable with respect to the collection of 
equilibrium solutions $x_1(t)\equiv\dots\equiv x_n(t)\equiv\rm 
constant$. 
\end{theo}

\begin{proof}
The proof of Theorem \ref{th1} is based on the evolution of the 
following set-valued function $\tV : X^{hn}\rightrightarrows X$,
\begin{equation}
\label{V}
\tV(\tx)\eqd \sigma(\pi(\tx))
\end{equation}
along the solutions of \eqref{equa}.
The fact that $t\mapsto\tV(\tx(t))$ is non-increasing is stated in the 
following result.
\begin{lemm}
\label{le2}
Let $x$ be a solution of equation \eqref{equa}.
Then, for all $t\in\Nset$,
\[
\tV\left(
\tx(t+1)
\right) \subseteq \tV\left(
\tx(t)
\right)\ .
\]
\end{lemm}
Let us first prove Lemma \ref{le2}.
For any $k\in\cN$, for any $t\in\Nset$,
\[
x_k(t+1) 
\in \sigma\left(
\{x_k(t)\}\cup\{x_i(t-j)\ :\ ((i,j),k)\in\cA(t)\}
\right) \subseteq \tV(\tx(t))
\]
successively by Assumption \ref{as1}
and Assumption \ref{as2}.\ref{assumption2},
and one concludes the demonstration of Lemma \ref{le2} by use of 
Assumptions \ref{as2}.\ref{assumption2} and 
\ref{as2}.\ref{assumption1}.
The proof of Theorem \ref{th1} is then obtained as a direct 
consequence.
\end{proof}

In view of Lemma \ref{le2}, one may now have a clearer understanding of 
the fact that the map $\sigma$ has a double role: it is necessary to
define the flow, but also serves as a set-valued Lyapunov function of 
the systems.
Indeed, Assumption \ref{as1} states that each agent has to remain in 
the set $\tV(\tx(t))$, of which it has only an imperfect knowledge, 
and does its best to come closer from the other agents it has detected 
(this is the meaning of the use of the relative interior).
In particular, when no new information is received, the only possible 
choice is to stay at the same place.

As detailed in Section \ref{se22}, contrary to $\sigma$, the map $\ri\ 
\sigma$ is not monotone:
violation of this rule may occur when $S' \subset S$ 
and the $\sigma$-hulls $\sigma(S), \sigma(S')$ have different 
topological dimensions as spheres.
Up to this subtlety, a consequence of Assumption \ref{as1} is that,  in 
general, {\em the larger the 
quantity of information received by agent $k$ from its neighborhood, 
the largest the set of possible updates it may choose\/}
(see  the monotony property in Assumption \ref{as2}.\ref{assumption2}).
Although this may sound paradoxical at first glance, this increase of the 
decision possibilities is quite natural: it means that supplementary 
information either leads to make a choice which could have been done 
otherwise (it is ignored or makes more valuable the decision) or allows 
to 
adopt choices which would not have been done otherwise.
The ``subtlety'' comes from the fact that, when the information 
available to an agent is poor, some decisions are taken which would not
 have been possible with richer data.
For example, the possibility of staying in the same place, which occurs 
when an agent, say agent $1$, is isolated from the other world,
disappears when the position of another agent located elsewhere, agent 
$2$, is received.
However, the unique choice $\sigma(\{x_1\})=\{x_1\}$ is then located 
``on the boundary'' of the decision set $\ri\ \sigma(\{x_1,x_2\})$,
 see Lemma \ref{le17}.\\

The key result of the paper is now stated.
It provides a {\em necessary and sufficient stability condition\/} for 
system \eqref{equa}, which extends \cite[Theorem 3]{MOR03}.

\begin{theo}
\label{th2}
Assume that Assumptions \ref{as1} and \ref{as2} are fulfilled.
Then the discrete-time system \eqref{equa} is uniformly globally 
attractive with respect to the collection of equilibrium solutions 
$x_1(t)\equiv\dots\equiv x_n(t)\equiv\rm constant$ {\em if and only 
if\/} there 
exists
$T\geq 0$ such that for all $t_0\in\Nset$ there is a node connected to 
all other nodes across $[t_0,t_0+T]$.
\end{theo}

The uniformity which is meant in the statement of Theorems \ref{th1} 
and \ref{th2} is with respect to {\em time}.
One may check from the proofs that it is also valid with respect to the 
different trajectories of \eqref{equa}.

\begin{proof}[Proof of Theorem \ref{th2}]
(Only if.)
The proof consists in an adaptation of the contraposition argument 
developed by Moreau \cite[Proof of Theorem 3]{MOR03}.
Assume that for every $T\geq 0$ there is $t_0\in\Nset$ such
that the 
sequence of graphs with delays has no node connected to each other 
across 
the interval $[t_0, t_0 + T]$.
This implies that for every $T\geq 0$ there is $t_0\in\Nset$ and 
nonempty, disjoint
subsets $\cL_1, \cL_2\subset\cN$ such that $\neig(\cL_1, 
\cA(t))=\neig(\cL_2,\cA(t))=\emptyset$ for all $t\in [t_0, t_0 + T]$.
The proof of this fact consists in checking that the proof of 
\cite[Theorem 5]{MOR03} holds also in the case of graph with delays as 
defined 
above in Definition \ref{de1}.

Let $y, \by\in X$ and consider any solution $\zeta$ of \eqref{equa} 
departing from initial data defined by:
\[
\zeta_k(t_0-j)\begin{cases}
= y &\forall (k,j)\in \cL_1\times\cH,\\
= \by &\forall (k,j)\in \cL_2\times\cH,\\
\in\sigma(\{ y, \by\}) &\forall (k,j)\in (\cN\setminus 
(\cL_1\cup\cL_2)) 
\times\cH \ .
\end{cases}
\]
As in the proof of \cite[Theorem 5]{MOR03}, we still have the same 
relation at time $t_0 + T + 1$, since $\neig(\cL_1, 
\cA(t))=\neig(\cL_2,\cA(t))=\emptyset$ for all $t\in [t_0, t_0 + T]$.
As the time $T$ may be chosen arbitrarily large, this contradicts 
uniform global attractivity
of \eqref{equa} with respect to the equilibrium 
solutions $x_1(t)\equiv\dots\equiv x_n(t)\equiv\textrm{constant}$.

(If.)
As in \cite[Theorem 3]{MOR03}, the proof is based on a (strict) 
decrease property of the set-valued function $\tV$ introduced in 
\eqref{V}. 

Let $T\geq 0$ chosen as in the statement of Theorem \ref{th2}, $x$ an 
arbitrary solution of \eqref{equa}, and $t_0\in\Nset$ for which the 
values of $x_k(t-j)$, $k\in\cN$, $j\in\cH$ are not all equal.
For all $k \in \cN$, define the integer-valued function $\alpha_k (t)$, 
$t\geq t_0$, by:

\[
\alpha_{k}(t) 
\eqd\card\ \{ j \in\cN :\ x_{j}(t)\in\Sigma_{x_k (t_0) } \}\ .
\]
where $\Sigma_x$ is meant relative to the set $S=\pi ( \tx (t_0) )$.
In 
words, this is the number of
agents which at time $t \geq t_0$ are still belonging to the critical 
portion of the boundary $\Sigma_{x_k(t_0)}$.
Assume by contradiction that at time $t_1>t_0$ a new agent $x_\alpha$ 
enters $\Sigma_{x_k(t_0)}$ which was not there at time $t_1-1$
($x_\alpha (t_1 -1) \notin \Sigma_{x_k (t_0)}$).
Let $S'$ denote the set of points in $X$
used by $x_\alpha$ at time $t_1-1$ in order to
update its state.
Of course, by Assumption \ref{as1}, $S'$ comprises 
$x_\alpha (t_1-1)$ itself.
Moreover, by monotonicity of the set-valued Lyapunov function $\tV$, 
$S' \subset \sigma (S)$. Now, by the updating rule \ref{as1},
$x_\alpha (t_1) \in \Sigma_{x_k (t_0)}$ is only possible provided that 
$\ri ( \sigma(S') ) \cap \Sigma_{x_k (t_0)} \neq
\emptyset$, and therefore, application of Assumption 
\ref{as2}.\ref{assumption3a}
yields 
$x_\alpha (t_1-1) \in S' \subset  \Sigma_{x_k (t_0)}$,
which contradicts what was
previously stated.
Hence, we can conclude that 
agents can only leave $\Sigma_{x_k (t_0)}$, but never get back in.  
In particular then, the functions $\alpha_k$ satisfy the inequality:
\[
\forall t\in [t_0,+\infty),\forall \, k \in\cN,\
\alpha_{k}(t+1)\leq \alpha_{k}(t)
\]
and
\begin{equation}
\label{tobequotedlater}
\alpha_k (t_0) = \alpha_k (t_1) 
\Leftrightarrow \{ j \in 
\cN: \, x_j (t_0) \in\Sigma_{x_k (t_0)} \} =
 \{ j \in \cN: \, x_j (t) \in\Sigma_{x_k(t_0)} \} \quad \forall \, t \in 
\{ t_0,\ldots , t_1 \}. 
\end{equation}
On the other hand, if $t\geq t_0$ is such that $\alpha_k(t)=0$ for a 
certain $k$ in $\cN$, then Assumption \ref{as2}.\ref{assumption3b}
 implies that $\mu(\sigma(\tx(t)))<\mu(\sigma(\tx(t_0)))$.
 An important 
step consists in showing that:
 \begin{equation}
 \label{decr}t> t_0 + T' \Rightarrow \exists k\in\cN, 
\alpha_k(t)<\alpha_k(t_0),\quad T'\eqd h+T\ .
 \end{equation}
 There are at most $n$ different sets $\Sigma_k$ in $\sigma(\tx(t_0))$, 
and $\alpha_k(t_0)\leq n-1$.
 Consequently, the repetition of the argument used to get implication 
\eqref{decr} (if allowed) will yield:
 \[ t> t_0 + (n-1)^2T' \Rightarrow\exists k\in\cN, \alpha_k(t)=0\ .\]
 As a consequence, the estimate:
 \begin{equation}
 \label{yiel}t> t_0 + T'' \Rightarrow 
\mu(\tV(\tx(t)))<\mu(\tV(\tx(t_0))),\quad T''\eqd (n-1)^2T'\ ,
 \end{equation}
 will be deduced from Assumption \ref{as2}.\ref{assumption3b}, because $\pi(\tx(t))$ being a {\em finite\/} set of points located in $\sigma(\pi(\tx(t_0)))\setminus\Sigma_{x_k(t_0)}$, is thus at a {\em nonzero\/} distance (more precisely a $\sigma(\pi(\tx(t_0)))$-distance, see \eqref{last}) from $\Sigma_{x_k(t_0)}$.
 In order 
to get \eqref{yiel}, let us now prove \eqref{decr}.
 \begin{enumerate}
 \item Inequality \eqref{decr} is fulfilled if it holds for $t < 
t_0+h+T$.
 \item Otherwise, one has $\alpha_k (t_0+h+T)= \alpha_k(t_0)$ for all 
$k$ in $\cN$ and,
 by virtue of (\ref{tobequotedlater}), the set $\cL_k \eqd \{ j \in 
\cN: x_j (t) \in \Sigma_{x_k (t_0)} \}$ has not changed
 for $t \in \{t_0, \ldots , t_0 + h +T \}$.

 Using the hypothesis in the statement of Theorem \ref{th2}, there 
exists an agent, numbered $k$, connected to all others
 across the interval $[t_0+h, t_0+h+T]$. 
  By definition, the set $\Sigma_{x_k (t_0) }$ does not contain $x_k (t_0)$,
  and,
since $\cL_k$ has not varied in time, then also $x_k(t)\not\in
  \Sigma_{x_k(t_0)}$ for 
 $t=t_0, \dots, t_0+h+T$.
 Moreover, because of the connectivity property of the 
graph put in the statement,
  $\neig(\cL_k,\cup_{t\in [t_0+h,t_0+h+T]}\cA(t))\neq\emptyset$.
  Let $i \in\cL_k$ be such that $\neig( i ,\cup_{t\in 
[t_0+h,t_0+h+T]}\cA(t))\setminus\cL_k\neq\emptyset$,
viz.\ an agent which over the time-interval $[t_0+h,t_0+h+T]$ is receiving 
information from outside $\cL_k$.
  Clearly $x_i(t_0+h-1)\in\Sigma_{x_k (t_0) }$, but 
$x_i(t_0+h+T)\not\in\Sigma_{x_k (t_0)}$, as Assumption 
  \ref{as2}.\ref{assumption3a} implies that elements in $\Sigma_x$ are 
attracted outside of it, as soon as they
  receive information from agents sitting outside $\Sigma_x$.
   This yields finally: 
   $\alpha_k(t_0+h+T)< \alpha_k(t_0)$ for the value of $k$ previously 
exhibited. Inequality \eqref{decr} is thus proved.
   Of course this is true only as long as $\Sigma_{x_k (t_0)} $ is 
non-empty to start with. 
   \end{enumerate}

One verifies easily that the same argument may be used recursively, 
because the sets $\Sigma_{x_k (t_0)}$ may be kept unchanged as long as 
$\alpha_k(t)>0$ for all $k\in\cN$.
Thus, \eqref{yiel} is proved.\\

Considering now $\tx(t_0)\in X^{hn}$ as a variable, let
\[
\beta(\tx(t_0))\eqd\inf
\mu\left( \tV(\zeta(0))
\right)
- \mu\left(
\tV(\zeta(T''))
\right)\ ,
\]
where the infimum is taken over all sequences $\zeta(1),\dots, 
\zeta(T'')$ in $X^{hn}$ such that $\zeta(0)=\tx(t_0)$ 
and, for 
all $ t=1,\dots, T''$, for all $k\in\cN$,
\[
\zeta_{k,0}(t)\in e_k(\cA(t_0+t))(\zeta(t-1))
 \text{ and }
\zeta_{k,j}(t)= \zeta_{k,j-1}(t-1) \text{ for } j\in\cH\setminus\{0\}\ 
.
\]

The meaning of the previous line is precisely that the infimum is 
computed over {\em all possible trajectories\/}
of the difference inclusion \eqref{equa}.
Now, the collection of $\zeta(t)\in X^{hn}$, $t=1,\dots, T''$ 
satisfying the previous condition is nonempty and compact for all 
initial value $\tx(t_0)\in X^{hn}$.
Indeed,
by Assumption \ref{as1}, the set-valued functions $e_k(\cA)$ are continuous 
and take nonempty, compact values.

The quantity to be minimized is strictly positive when the $hn$ 
components of $\zeta(0)$ are not all equal, due to the strict decrease 
property of $\tV$ established above and Assumption 
\ref{as2}.\ref{assumption3}.
Also, the expression to be minimized is lower semicontinuous with 
respect to $\zeta(1),\dots, \zeta(T'')$, as $\tV$ is 
continuous
(by Assumption \ref{as2}) and $\mu$ is lower semicontinuous (by Lemma 
\ref{le15}).
Thus,  $\beta(\tx(t_0))>0$, except if all the components of 
$\tx(t_0)$ are all equal.
In other words, $\beta$ is definite positive with respect to $\{\tx\in 
X^{hn}\ :\ x_1=\dots=x_{hn}\}$.

By Assumption \ref{as2}.\ref{assumption4}, the map
$X^{hn}\to\Rset^+$, $\tx\mapsto\mu(\tV(\tx))$ is continuous.
The proof of Theorem \ref{th2} is then achieved as for \cite[Theorem 
3]{MOR03},
by use of a result on set-valued 
Lyapunov functions.
The latter, Theorem \ref{th3}, is an extension of \cite[Theorem 1]{MOR03} to differential inclusions, given in Appendix.
\end{proof}

\begin{exam}
\label{ex4}
The necessity for each agent to take into account the present values of 
its own position may be seen by the following counter-example.
We take $n=3$ and $h=2$.
Let
\[
\cA(2t) =( ((2,1),1),((1,0),2)),\
\cA(2t+1) =( ((2,1),3),((3,0),2)) \ .
\]
In other words, agent 2 sends alternativaly to agent 1 and 3 the value 
of its position at the previous instant, and receives in the same time 
the present value of their position, see Figure \ref{fi8}.
 \begin{figure}
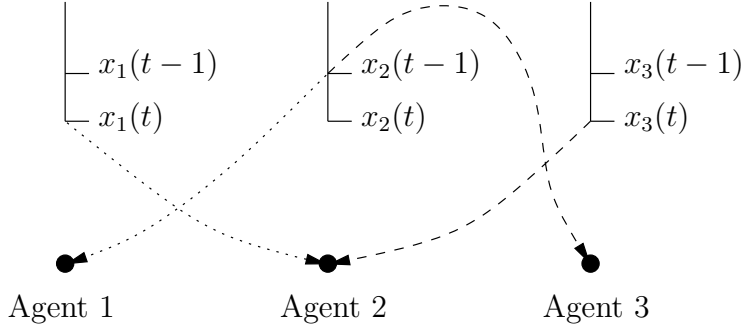

 \begin{center}
\input fig8.pstex_t
  \end{center} 
 \caption{Graph with delays representing the information flow for 
Example 
\ref{ex4}: even (dots) and odd (dash) times.}
   \label{fi8}
 \end{figure}
 Assume the agents use at time $t$ the value of their position at time 
$t-1$ to elaborate the update applied at time $t+1$.
Clearly, for the corresponding graph with delays, the agent 2 is 
connected to all other agents across 
any interval $[t,t+1]$.
However, one sees easily that provided that the agents 1 and 3 are 
located initially at different positions, the positions of agent 2 at 
even 
and odd times tend in general toward two different values.
As indicated by the existence of periodic motions, the strict decrease 
of the map $t\mapsto \mu(\tV(\tx(t)))$ may fail.
\end{exam}

\appendix

\addtocounter{section}{1}
\section*{Appendix \thesection\ -- Proof of Lemma \ref{le15}}

1.\
We first recall that the length of a Lipschitz arc $\psi(\lambda)$ 
defined on $[0,1]$ is equal to
\[
\leng(\psi) \eqd\int_0^1 \left\|
\frac{d\psi}{d\lambda}
\right\|\cdot d\lambda\ .\]
Notice that,
by Rademacher's theorem, $\psi$ Lipschitz implies differentiability 
almost everywhere, and therefore the
previous integral is well defined.
Let $x^0,x^1\in S$.
Taking $\varphi : X\to X$ as in \eqref{cS}, define the map $\psi : 
[0,1]\to X$
\begin{equation}
\label{psi}
\psi(\lambda)\eqd \varphi^{-1}\left(
(1-\lambda)\varphi(x^0)+\lambda\varphi(x^1)
\right)\ .
\end{equation}
As $\varphi(S)$ is convex, $\psi$ maps $[0,1]$ in $S$.
Moreover, due to the regularity assumption on $\varphi$, it is a 
Lipschitz 
arc, and
$ \psi(0)=x^0, \psi(1)=x^1$.
Thus the set
$\left\{
\leng(\psi)\ :\ \psi: 
[0,1]\stackrel{\textrm{Lipschitz}}{\longrightarrow}S, 
\psi(0)=x^0, \psi(1)=x^1
\right\}$
is non-void, and the definition of the map $d_S(x^0,x^1)$ given in the 
statement is meaningful.

Let us show its continuity with respect to $x^0,x^1\in S$.
Let $(x^0,x^1)\in S\times S$.
Consider a sequence $(x^\varepsilon)_{\varepsilon>0}$ of elements of 
$S$ tending towards $x^0$.
Let $\psi^{\varepsilon| 0}$ be a fixed Lipschitz arc linking 
$x^\varepsilon$ to $x^0$.
For any piecewise Lipschitz arc $\psi^{0| 1}$ linking $x^0$ to $x^1$, 
one 
may construct by concatenation of 
$\psi^{\varepsilon| 0}$ and $\psi^{0| 1}$ another Lipschitz arc 
$\psi^{\varepsilon| 1}$ linking $x^\varepsilon$ to $x^1$.
One has
\[
\leng(\psi^{\varepsilon| 1})
= \leng(\psi^{\varepsilon| 0})+\leng(\psi^{0|1})\ ,
\]
so
\[
d_S (x^{\varepsilon},x^1) \leq d_S (x^{\varepsilon},x^0)+ d_S 
(x^0,x^1)\ .
\]
Arguing similarly, one gets that $| d_S (x^{\varepsilon}, x^1 ) - d_S ( 
x^0,x^1) | 
\leq d_S (x^\varepsilon,x^0)$.

On the other hand, one may take $\psi$ as in \eqref{psi}, in such a way 
that
\[
\inf\left\{
\leng(\psi)\ :\ \psi: 
[0,1]\stackrel{\textrm{Lipschitz}}{\longrightarrow}S, 
\psi(0)=x^\varepsilon, 
\psi(1)=x^0
\right\}
\leq \left\|
\frac{d\psi}{d\lambda}
\right\|_{L^\infty} \|x^\varepsilon-x^0\|\ ,
\]
which shows the desired continuity property.

Defining $\mu(S)$ as in \eqref{mu}, one has $\mu(S)<+\infty$, because 
the image of a compact set by a continuous function is bounded.
Moreover, if $\mu(S)=0$, then, for any $(x^0,x^1)\in S\times S$, the 
length defined
 by the map $\psi$ in \eqref{psi} is zero, that is $x^0=x^1$ and $S$ is a 
singleton.
Conversely, if  $S$ is a singleton, then $\mu(S)=0$.
Last, for any Lipschitz arc $\psi$ linking $x^0$ to $x^1$ and defined 
as in 
\eqref{psi},
\[
\leng(\psi)
= \int_0^1 \left\|
\frac{d\psi}{d\lambda}
\right\|\cdot d\lambda
\geq \left\|
 \int_0^1\frac{d\psi}{d\lambda}\cdot d\lambda
\right\|
= \|x^0-x^1\|\ ,
\]
and this shows that $\mu(S)$ is at least equal to the maximal euclidian 
distance between two points of $S$, that is its diameter.
The equality when $S$ is convex is straightforward, taking the identity 
for $\varphi$ in \eqref{psi}.

We now prove the lower semicontinuity of $\mu$.
Let $S\in\cS$, and a sequence of sets $S_n\in\cS$ tending towards $S$ 
for the topology induced by Hausdorff distance.
Our goal is to prove that:
\begin{equation}
\label{geq}
\liminf_{n\to +\infty} \mu(S_n) \geq \mu(S)\ .
\end{equation}
Obviously, in order to establish inequality \eqref{geq}, it is 
sufficient to consider only sets $S_n$ such that $\mu(S_n)$ is bounded 
from above by a given constant, say by twice the value of $\mu(S)$.
Let $\varepsilon >0$, and consider two arbitrary sequences $x^0_n, 
x^1_n\in S_n$ and a sequence 
of Lipschitz arcs $\psi_n$ linking $x^0_n$ to $x^1_n$ 
in $S_n$ and of length at most equal to 
$d_{S_n}(x^0_n,x^1_n)+\varepsilon$.
We thus have:
\begin{equation}
\label{inea}
\liminf_{n \rightarrow + \infty} d_{S_n}(x^0_n,x^1_n)+\varepsilon \geq 
\liminf_{n \rightarrow + 
\infty}  \leng (\psi_n)
\geq \liminf_{n \rightarrow + \infty} d_{S_n}(x^0_n,x^1_n)\ .
\end{equation}
However, due to the previous remark on the boundedness of the sequence 
$\mu(S_n)$, one may assume without loss of generality that the arcs 
$\psi_n$ are 
covered with a rate of variation 
$\|\frac{d\psi_n}{d\lambda}\|_{L^\infty}$ 
uniformly bounded.
Indeed, if this is not the case, replace $\psi_n$ by the map $\tpsi_n$ 
defined by
\[
\tpsi_n \eqd \psi_n\circ \theta^{-1},\
\theta (t) \eqd \frac{\int_0^t \left\|
\frac{d\psi_n}{d\lambda}
\right\|\cdot d\lambda}
{\int_0^1 \left\|
\frac{d\psi_n}{d\lambda}
\right\|\cdot d\lambda}\ .
\]
The map $\tpsi_n$ has the same image and length than $\psi_n$, but the 
norm of its derivative is equal almost everywhere on $[0,1]$ to the 
constant $\int_0^1 \left\|
\frac{d\psi_n}{d\lambda}
\right\|\cdot d\lambda$.
In particular,
\begin{equation}
\label{ineb}
\liminf_{n \rightarrow + 
\infty}  \leng (\tpsi_n)=\liminf_{n \rightarrow + 
\infty}  \leng (\psi_n)\ .
\end{equation}
The arcs considered at this stage are thus equicontinuous.
By compactness, one deduces that there exist subsequences (denoted 
similarly $x^0_n, x^1_n$ and $\tpsi_n$) such that
\[
x^0_n\to x^0\in S,\ x^1_n\to x^1\in S,\ \tpsi_n\to\tpsi\in 
\textrm{Lipschitz} ([0,1];S)\ .
\]
In particular, since by Arzela-Ascoli $\tpsi_n \to \tpsi$ uniformly 
over compact sets, we also have
$\leng ( \tpsi_n ) \to \leng ( \tpsi )$. 
Thus, since $\mu (S_n) \geq d_{S_n} ( x^0_n , x^1_n)$, we have:
\begin{eqnarray}
\label{inequalitychain}
 \liminf_{n \rightarrow + \infty} \mu (S_n) + \varepsilon &\geq&
\lim_{n \rightarrow + 
\infty} d_{S_n} ( x^0_n,x^1_n ) + \varepsilon \nonumber \\ &\geq&
 \liminf_{n\rightarrow + \infty} \leng( \tpsi_n )
= \leng ( \tpsi ) \nonumber \\ &\geq &d_S ( x^0,x^1 ) 
\end{eqnarray}
By arbitrariety of $\varepsilon$:
\[ \liminf_{n \rightarrow + \infty} \mu (S_n) \geq d_S ( x^0,x^1 ). \]
We finally use arbitrariety of $x^0$,$x^1$ in $S$ (arbitrary converging 
sequences in $S_n$ yield arbitrary limit points in $S$ since
$S_n \rightarrow S$) to conclude $\liminf_{n \rightarrow + \infty} \mu 
(S_n) 
\geq \mu (S)$.
The lower semicontinuity of $\mu$ is demonstrated.\\

\noindent 2.\
To prove that the definition of $\ri(S)$ is independent of $\varphi$ 
amounts to show that: if the image of a convex set $S$ by a bijective 
bi-Lipschitz map $f: X\to X$ is a convex set, then the image of the 
relative 
interior $\ri(S)$ is the relative interior of the image of $S$.
Let $x\in \ri(S)$.
Consider the restriction of $f$ to the affine hull $\ah(S)$ of $S$.
The set $\ri(S)$ is convex, thus there exists a convex neighborhood 
$V\subset S$ of $x$ in $\ah(S)$.
By continuity of $f^{-1}$, the image of $V$ through $f$ is a 
neighborhood of $f(x)$ in $f(\ah(S))$.
By hypothesis, $f(S)$ is convex, thus there exists a convex 
neighborhood $W\subset f(S)$ of $f(x)$ in $\ah(f(S))$.
Now, the intersection $f(V)\cap W\subset f(S)$ is a neighborhood of 
$f(x)$ in $\ah(f(S))$, so $f(x)\in\ri(f(S))$.
We thus get $\ri(S)\subseteq f^{-1}(\ri(f(S)))$, and one shows 
similarly the converse inclusion.

Let $S\in\cS$.
If $\ri(S)=\emptyset$, then $\ri(\varphi(S))=\emptyset$ and $S$ is a 
singleton.
Conversely, if $S$ is a singleton, $\varphi(S)$ is also a singleton and 
$\ri(\varphi(S))=\emptyset$.

It is clear that $\ri\ (S)\subset S$.
Also, the fact that $\inte\ (\varphi(S))\subseteq\ri\ 
(\varphi(S))\subset \varphi(S)$ implies that
\begin{eqnarray*}
\inte\ S
&= &
\inte\ \varphi^{-1}\left( \varphi(S) 
\right)\subseteq\varphi^{-1}(\inte\ (\varphi(S)))
\quad \text{(as $\varphi^{-1}$ is continuous)}\\
&\subseteq &
\varphi^{-1}(\ri\ (\varphi(S)))
=
\ri\ (S)
\quad\text{ (by definition of $\ri\ (S)$)}\ .
\end{eqnarray*}

Last, when $S$ is a convex set, one may take for $\varphi$ in 
\eqref{cS} the identity, and this proves that in this case $\ri(S)$ is 
equal to 
the relative interior of $S$.

This ends the proof of Lemma \ref{le15}.

\addtocounter{section}{1}
\section*{Appendix \thesection\ -- Stability based on set-valued 
Lyapunov functions}

The following result is an adaptation of \cite[Theorem 1]{MOR03} to difference inclusions.
For sake of completeness, a proof is provided, intimately linked to the proof of Moreau's result.

\begin{theo}
\label{th3}
Let $\cX$ be a finite-dimensional Euclidean space and consider a 
continuous set-valued map $e\ :\ \Nset\times\cX\rightrightarrows \cX$ taking on closed values,  giving rise to the difference inclusion \eqref{eque}.
Let $\Xi$ be a collection of equilibrium solutions and denote the 
corresponding set of equilibrium points by $\Phi$.
Consider an upper semicontinuous \cite[p.\ 41]{aubin}
set-valued function $V : X 
\rightrightarrows X$ satisfying
\begin{enumerate}
\item
$x \in V (x), \forall x \in\cX$;
\item
$\displaystyle\bigcup_{y\in e(t,x)} V (y) \subseteq V (x), \forall t\in\Nset, \forall x\in\cX$.
\end{enumerate}

If $V (\phi) = \{\phi\}$ for all $\phi\in\Phi$ then the dynamical 
system  is uniformly stable with respect to $\Phi$.
If $V (x)$ is bounded for all $x\in\cX$ then the dynamical system is 
uniformly bounded with respect to $\Phi$.

Consider in addition a function 
$\mu : \textrm{Image}(V ) \to\Rset_{\geq 0}$ and a lower semicontinuous 
function $\beta : \cX \to\Rset_{\geq 0}$ satisfying
\begin{enumerate}
\item[3.] $\mu\circ V : \cX\to\Rset_{\geq 0}$ is bounded uniformly with 
respect to $x$ in bounded subsets of $\cX$;
\item[4.] $\beta$ is 
positive definite with respect to $\Phi$ that is, $\beta(\phi) = 0$ for all 
$\phi\in\Phi$ and $\beta(x) > 0$ for all $x\in\cX\setminus\Phi$;
\item[5.] 
$\displaystyle\sup_{y\in e(t,x)}\mu(V (y)) -\mu(V (x)) \leq -\beta(x), \forall t\in\Nset, 
\forall x\in\cX$.
\end{enumerate}
If $V (\phi) = \{\phi\}$ for all $\phi\in\Phi$ and $V (x)$ is bounded 
for all $x\in\cX$ then the dynamical system is uniformly globally 
asymptotically stable with respect to $\Phi$.
\end{theo}

The results stated above remain true if, for a fixed $\tau\in\Nset$,
 the decrease relations in 2.\ and 5.\ occur between $V(y_{i+1})$ and $V(y_i)$, $i=0, \dots, \tau-1$,
 $y_0=x$, $y_\tau=y$, instead of $V(y)$ and $V(x)$.

\begin{proof}
(Uniform stability.)  Consider arbitrary $\varphi\in\Phi$
and $c_2 > 0$.
If $V (\varphi) = \{\varphi\}$ then, by upper semicontinuity
of $V$, there is $c_1 > 0$ such that $V (x) \subset B(\varphi, c_2)$ for all $x\in B(\varphi, c_1)$.
Consider arbitrary $t_0 \in\Nset$ and $x_0 \in B(\varphi, c_1)$
and let $\zeta$ denote any solution of inclusion  \eqref{eque} with $\zeta(t_0) = x_0$.
Conditions 1 and 2 of Theorem \ref{th3} imply that, for all $t\geq t_0$,
\[
\zeta(t) \in \bigcup_{y\in e(t,x_0)} V (y) \subseteq V (x_0)\subset B(\varphi, c_2)\ .
\]

(Uniform boundedness.) Consider arbitrary $\varphi\in\Phi$
and $c_1 > 0$.
If $V (x)$ is bounded for all $x \in \cX$ then, by upper semicontinuity of $V$, there is $c_2 > 0$ such that
$V (x) \subset B(\varphi, c_2)$ for all $x \in B(\varphi, c_1)$.
Consider arbitrary $t_0 \in \Nset$ and $x_0 \in B(\varphi, c_1)$ and let $\zeta$ be any solution of \eqref{eque}
with $\zeta(t_0) = x_0$.
Conditions 1 and 2 of Theorem \ref{th3}  imply
that for all $t\geq t_0$,
\[
\zeta(t) \in \bigcup_{y\in e(t,x_0)} V (y) \subseteq V (x_0)\subset B(\varphi, c_2)\ .
\]

(Uniform global asymptotic stability.)
It remains only to prove uniform global attractivity with respect
to $\Xi$.

Consider arbitrary $\varphi_1 \in \Phi$ and $c_1 > 0$.
If $V (x)$ is
bounded for all $x \in \cX$ then, by upper semicontinuity of $V$,
there is a compact set $K \subset\cX$ such that $V (x) \subseteq K$ for
all $x \in B(\varphi_1, c_1)$.
Similarly as above, Conditions 1 and 2 of
Theorem \ref{th3}  imply that every solution of \eqref{eque} initiated
in $B(\varphi_1, c_1)$ remains in $K$.

Consider in addition arbitrary $c_2 > 0$.
If $V (\varphi) =\{\varphi\}$ for
all $\varphi \in \Phi$ then, by upper semicontinuity of $V$, there is $c_3 >
0$ such that for all $x \in B(\Phi \cap K, c_3)$ there is $\varphi_2 \in \Phi$ such
that $V (x) \subset B(\varphi_2, c_2)$.
Similarly as above, Conditions 1
and 2 of Theorem \ref{th3}  imply that every solution of \eqref{eque}
entering $B(\Phi \cap K, c_3)$ remains in a $c_2$-ball around some
equilibrium point $\varphi_2 \in \Phi$.

It remains to prove the existence of $T \geq 0$ such that
every solution of \eqref{eque} starting in $B(\varphi_1, c_1)$ cannot remain
longer than $T$ subsequent times in $K$ without entering
$B(\Phi \cap K, c_3)$.
In agreement with Conditions 3 and 4 of
Theorem \ref{th3}  and the lower semicontinuity of $\beta$, we introduce
two real numbers:
\[
M \eqd \sup_{x\in B(\varphi_1,c_1)} \mu(V (x)) < \infty
\quad\text{ and }\quad
\Delta \eqd \min_{x \in K\setminus B(\Phi,c_3)} \beta(x) > 0\ .
\]

Let $T \geq 0$ be such that $T\Delta > M$.
Consider arbitrary
$t_0 \in \Nset$ and $x_0 \in B(\varphi_1, c_1)$ and let $\zeta$ denote any
solution of \eqref{eque} with $\zeta(t_0) = x_0$.
Then Condition 5
of Theorem \ref{th3}  implies that for some $t_1 \in [t_0, t_0 + T]$,
$\zeta(t_1) \in B(\Phi \cap K, c_3)$,
since otherwise $\zeta(t)\in K\setminus B(\Phi,c_3)$ for all $t \in [t_0, t_0 + T]$ and
\[
\mu(V (\zeta(t_0 + T)))\leq \mu(V(\zeta(t_0))) -T\min_{x \in K\setminus B(\Phi,c_3)} \beta(x)
\leq M-T\Delta < 0\ ,
\]
contradicting that $\mu$ takes only non-negative values.
Putting everything together, we conclude that for
some $\varphi_2 \in \Phi$ and for all $t\geq t_0+T$,
\[
\zeta(t) \in V (\zeta(t)) \subseteq \mu(V (\zeta(t_1)) \subset B(\varphi_2, c_2)\ .
\]
This achieves the proof of Theorem \ref{th3}.
\end{proof}

\end{document}

%% file: fig7.pstex_t
\begin{picture}(0,0)%
\epsfig{file=fig7.pstex}%
\end{picture}%
\setlength{\unitlength}{3947sp}%
\begingroup\makeatletter\ifx\SetFigFont\undefined%
\gdef\SetFigFont#1#2#3#4#5{%
  \reset@font\fontsize{#1}{#2pt}%
  \fontfamily{#3}\fontseries{#4}\fontshape{#5}%
  \selectfont}%
\fi\endgroup%
\begin{picture}(4804,3316)(4037,-5465)
\put(6361,-2461){\makebox(0,0)[lb]{\smash{{\SetFigFont{12}{14.4}{\rmdefault}{\mddefault}{\updefault}$x_2(t-5)$}}}}
\put(6361,-3061){\makebox(0,0)[lb]{\smash{{\SetFigFont{12}{14.4}{\rmdefault}{\mddefault}{\updefault}$x_2(t-3)$}}}}
\put(6361,-3361){\makebox(0,0)[lb]{\smash{{\SetFigFont{12}{14.4}{\rmdefault}{\mddefault}{\updefault}$x_2(t-2)$}}}}
\put(6361,-3661){\makebox(0,0)[lb]{\smash{{\SetFigFont{12}{14.4}{\rmdefault}{\mddefault}{\updefault}$x_2(t-1)$}}}}
\put(6361,-3991){\makebox(0,0)[lb]{\smash{{\SetFigFont{12}{14.4}{\rmdefault}{\mddefault}{\updefault}$x_2(t)$}}}}
\put(8011,-3991){\makebox(0,0)[lb]{\smash{{\SetFigFont{12}{14.4}{\rmdefault}{\mddefault}{\updefault}$x_3(t)$}}}}
\put(8011,-2461){\makebox(0,0)[lb]{\smash{{\SetFigFont{12}{14.4}{\rmdefault}{\mddefault}{\updefault}$x_3(t-5)$}}}}
\put(8011,-3361){\makebox(0,0)[lb]{\smash{{\SetFigFont{12}{14.4}{\rmdefault}{\mddefault}{\updefault}$x_3(t-2)$}}}}
\put(8011,-3661){\makebox(0,0)[lb]{\smash{{\SetFigFont{12}{14.4}{\rmdefault}{\mddefault}{\updefault}$x_3(t-1)$}}}}
\put(8011,-2761){\makebox(0,0)[lb]{\smash{{\SetFigFont{12}{14.4}{\rmdefault}{\mddefault}{\updefault}$x_3(t-4)$}}}}
\put(8011,-3061){\makebox(0,0)[lb]{\smash{{\SetFigFont{12}{14.4}{\rmdefault}{\mddefault}{\updefault}$x_3(t-3)$}}}}
\put(4711,-3991){\makebox(0,0)[lb]{\smash{{\SetFigFont{12}{14.4}{\rmdefault}{\mddefault}{\updefault}$x_1(t)$}}}}
\put(4711,-2761){\makebox(0,0)[lb]{\smash{{\SetFigFont{12}{14.4}{\rmdefault}{\mddefault}{\updefault}$x_1(t-4)$}}}}
\put(4711,-2461){\makebox(0,0)[lb]{\smash{{\SetFigFont{12}{14.4}{\rmdefault}{\mddefault}{\updefault}$x_1(t-5)$}}}}
\put(4711,-3061){\makebox(0,0)[lb]{\smash{{\SetFigFont{12}{14.4}{\rmdefault}{\mddefault}{\updefault}$x_1(t-3)$}}}}
\put(4711,-3361){\makebox(0,0)[lb]{\smash{{\SetFigFont{12}{14.4}{\rmdefault}{\mddefault}{\updefault}$x_1(t-2)$}}}}
\put(4711,-3661){\makebox(0,0)[lb]{\smash{{\SetFigFont{12}{14.4}{\rmdefault}{\mddefault}{\updefault}$x_1(t-1)$}}}}
\put(4186,-5411){\makebox(0,0)[lb]{\smash{{\SetFigFont{12}{14.4}{\rmdefault}{\mddefault}{\updefault}Agent 1}}}}
\put(5896,-5411){\makebox(0,0)[lb]{\smash{{\SetFigFont{12}{14.4}{\rmdefault}{\mddefault}{\updefault}Agent 2}}}}
\put(7546,-5411){\makebox(0,0)[lb]{\smash{{\SetFigFont{12}{14.4}{\rmdefault}{\mddefault}{\updefault}Agent 3}}}}
\put(6361,-2761){\makebox(0,0)[lb]{\smash{{\SetFigFont{12}{14.4}{\rmdefault}{\mddefault}{\updefault}$x_2(t-4)$}}}}
\end{picture}%

%% file: fig1.pstex_t
\begin{picture}(0,0)%
\epsfig{file=fig1.pstex}%
\end{picture}%
\setlength{\unitlength}{3947sp}%
\begingroup\makeatletter\ifx\SetFigFont\undefined%
\gdef\SetFigFont#1#2#3#4#5{%
  \reset@font\fontsize{#1}{#2pt}%
  \fontfamily{#3}\fontseries{#4}\fontshape{#5}%
  \selectfont}%
\fi\endgroup%
\begin{picture}(1961,3214)(1651,-4019)
\put(1651,-2461){\makebox(0,0)[lb]{\smash{{\SetFigFont{12}{14.4}{\rmdefault}{\mddefault}{\updefault}$x_3$}}}}
\put(3151,-961){\makebox(0,0)[lb]{\smash{{\SetFigFont{12}{14.4}{\rmdefault}{\mddefault}{\updefault}$x_2$}}}}
\put(3151,-3961){\makebox(0,0)[lb]{\smash{{\SetFigFont{12}{14.4}{\rmdefault}{\mddefault}{\updefault}$x_1$}}}}
\put(3076,-2386){\makebox(0,0)[lb]{\smash{{\SetFigFont{12}{14.4}{\rmdefault}{\mddefault}{\updefault}$\sigma(\tx)$}}}}
\end{picture}%

%% file: fig2.pstex_t
\begin{picture}(0,0)%
\epsfig{file=fig2.pstex}%
\end{picture}%
\setlength{\unitlength}{3947sp}%
\begingroup\makeatletter\ifx\SetFigFont\undefined%
\gdef\SetFigFont#1#2#3#4#5{%
  \reset@font\fontsize{#1}{#2pt}%
  \fontfamily{#3}\fontseries{#4}\fontshape{#5}%
  \selectfont}%
\fi\endgroup%
\begin{picture}(4223,3370)(139,-4019)
\put(3901,-1711){\makebox(0,0)[lb]{\smash{{\SetFigFont{12}{14.4}{\rmdefault}{\mddefault}{\updefault}$x_5$}}}}
\put(3151,-3961){\makebox(0,0)[lb]{\smash{{\SetFigFont{12}{14.4}{\rmdefault}{\mddefault}{\updefault}$x_1$}}}}
\put(3076,-2386){\makebox(0,0)[lb]{\smash{{\SetFigFont{12}{14.4}{\rmdefault}{\mddefault}{\updefault}$\sigma(\tx)$}}}}
\put(976,-3436){\makebox(0,0)[lb]{\smash{{\SetFigFont{12}{14.4}{\rmdefault}{\mddefault}{\updefault}$e_1$}}}}
\put(676,-2386){\makebox(0,0)[lb]{\smash{{\SetFigFont{12}{14.4}{\rmdefault}{\mddefault}{\updefault}$e_2$}}}}
\put(1651,-2461){\makebox(0,0)[lb]{\smash{{\SetFigFont{12}{14.4}{\rmdefault}{\mddefault}{\updefault}$x_3$}}}}
\put(3151,-826){\makebox(0,0)[lb]{\smash{{\SetFigFont{12}{14.4}{\rmdefault}{\mddefault}{\updefault}$x_2$}}}}
\put(3826,-1111){\makebox(0,0)[lb]{\smash{{\SetFigFont{12}{14.4}{\rmdefault}{\mddefault}{\updefault}$x_4$}}}}
\end{picture}%

%% file: fig3.pstex_t
\begin{picture}(0,0)%
\epsfig{file=fig3.pstex}%
\end{picture}%
\setlength{\unitlength}{3907sp}%
\begingroup\makeatletter\ifx\SetFigFont\undefined%
\gdef\SetFigFont#1#2#3#4#5{%
  \reset@font\fontsize{#1}{#2pt}%
  \fontfamily{#3}\fontseries{#4}\fontshape{#5}%
  \selectfont}%
\fi\endgroup%
\begin{picture}(7187,4677)(807,-4838)
\put(5926,-1411){\makebox(0,0)[lb]{\smash{{\SetFigFont{12}{14.4}{\rmdefault}{\mddefault}{\updefault}$x_4$}}}}
\put(1651,-2461){\makebox(0,0)[lb]{\smash{{\SetFigFont{12}{14.4}{\rmdefault}{\mddefault}{\updefault}$x_3$}}}}
\put(3151,-826){\makebox(0,0)[lb]{\smash{{\SetFigFont{12}{14.4}{\rmdefault}{\mddefault}{\updefault}$x_2$}}}}
\put(2866,-3781){\makebox(0,0)[lb]{\smash{{\SetFigFont{12}{14.4}{\rmdefault}{\mddefault}{\updefault}$x_1$}}}}
\put(4876,-2611){\makebox(0,0)[lb]{\smash{{\SetFigFont{12}{14.4}{\rmdefault}{\mddefault}{\updefault}$\sigma(\tx)$}}}}
\end{picture}%

%% file: fig5.pstex_t
\begin{picture}(0,0)%
\epsfig{file=fig5.pstex}%
\end{picture}%
\setlength{\unitlength}{3947sp}%
\begingroup\makeatletter\ifx\SetFigFont\undefined%
\gdef\SetFigFont#1#2#3#4#5{%
  \reset@font\fontsize{#1}{#2pt}%
  \fontfamily{#3}\fontseries{#4}\fontshape{#5}%
  \selectfont}%
\fi\endgroup%
\begin{picture}(9602,6786)(0,-5946)
\put(1466,-3336){\makebox(0,0)[lb]{\smash{{\SetFigFont{12}{14.4}{\rmdefault}{\mddefault}{\updefault}$x_3$}}}}
\put(6976,-2836){\makebox(0,0)[lb]{\smash{{\SetFigFont{12}{14.4}{\rmdefault}{\mddefault}{\updefault}$\sigma_\varphi(\tx)$}}}}
\put(4206, 89){\makebox(0,0)[lb]{\smash{{\SetFigFont{12}{14.4}{\rmdefault}{\mddefault}{\updefault}$x_2$}}}}
\put(5621,-4976){\makebox(0,0)[lb]{\smash{{\SetFigFont{12}{14.4}{\rmdefault}{\mddefault}{\updefault}$x_1$}}}}
\end{picture}%

%% file: fig6.pstex_t
\begin{picture}(0,0)%
\epsfig{file=fig6.pstex}%
\end{picture}%
\setlength{\unitlength}{3947sp}%
\begingroup\makeatletter\ifx\SetFigFont\undefined%
\gdef\SetFigFont#1#2#3#4#5{%
  \reset@font\fontsize{#1}{#2pt}%
  \fontfamily{#3}\fontseries{#4}\fontshape{#5}%
  \selectfont}%
\fi\endgroup%
\begin{picture}(6001,4201)(1248,-4838)
\put(3676,-2236){\makebox(0,0)[lb]{\smash{{\SetFigFont{12}{14.4}{\rmdefault}{\mddefault}{\updefault}$\sigma_\cap(\tx)$}}}}
\put(1651,-2461){\makebox(0,0)[lb]{\smash{{\SetFigFont{12}{14.4}{\rmdefault}{\mddefault}{\updefault}$x_3$}}}}
\put(3151,-826){\makebox(0,0)[lb]{\smash{{\SetFigFont{12}{14.4}{\rmdefault}{\mddefault}{\updefault}$x_2$}}}}
\put(2866,-3781){\makebox(0,0)[lb]{\smash{{\SetFigFont{12}{14.4}{\rmdefault}{\mddefault}{\updefault}$x_1$}}}}
\end{picture}%

%% file: fig8.pstex_t
\begin{picture}(0,0)%
\epsfig{file=fig8.pstex}%
\end{picture}%
\setlength{\unitlength}{3947sp}%
\begingroup\makeatletter\ifx\SetFigFont\undefined%
\gdef\SetFigFont#1#2#3#4#5{%
  \reset@font\fontsize{#1}{#2pt}%
  \fontfamily{#3}\fontseries{#4}\fontshape{#5}%
  \selectfont}%
\fi\endgroup%
\begin{picture}(4700,2046)(4141,-5245)
\put(6361,-3991){\makebox(0,0)[lb]{\smash{{\SetFigFont{12}{14.4}{\rmdefault}{\mddefault}{\updefault}$x_2(t)$}}}}
\put(8011,-3991){\makebox(0,0)[lb]{\smash{{\SetFigFont{12}{14.4}{\rmdefault}{\mddefault}{\updefault}$x_3(t)$}}}}
\put(8011,-3661){\makebox(0,0)[lb]{\smash{{\SetFigFont{12}{14.4}{\rmdefault}{\mddefault}{\updefault}$x_3(t-1)$}}}}
\put(4711,-3991){\makebox(0,0)[lb]{\smash{{\SetFigFont{12}{14.4}{\rmdefault}{\mddefault}{\updefault}$x_1(t)$}}}}
\put(4711,-3661){\makebox(0,0)[lb]{\smash{{\SetFigFont{12}{14.4}{\rmdefault}{\mddefault}{\updefault}$x_1(t-1)$}}}}
\put(4141,-5191){\makebox(0,0)[lb]{\smash{{\SetFigFont{12}{14.4}{\rmdefault}{\mddefault}{\updefault}Agent 1}}}}
\put(5851,-5191){\makebox(0,0)[lb]{\smash{{\SetFigFont{12}{14.4}{\rmdefault}{\mddefault}{\updefault}Agent 2}}}}
\put(7501,-5191){\makebox(0,0)[lb]{\smash{{\SetFigFont{12}{14.4}{\rmdefault}{\mddefault}{\updefault}Agent 3}}}}
\put(6361,-3661){\makebox(0,0)[lb]{\smash{{\SetFigFont{12}{14.4}{\rmdefault}{\mddefault}{\updefault}$x_2(t-1)$}}}}
\end{picture}%

%% file: attempt_sub.bbl
\begin{thebibliography}{99}

\bibitem{aubin} J.P.\ Aubin and A.\ Cellina (1984). \emph{Differential 
Inclusions:
 set-valued maps and viability theory},
Springer-Verlag, Berlin

\bibitem{FAXIFAC} J.A.\ Fax and R.M.\ Murray (2002). Information flow and 
cooperative
control of vehicle formations.
In \emph{Proc.\ of the IFAC World 
Congress 2002},
Barcelona, Spain

\bibitem{Hirsch} M.W.\ Hirsch (1989). Convergent activation dynamics in
continuous time networks.
In \emph{Neural Networks}, Vol.\ 2, pp.\ 
331--349

\bibitem{Jad2003} A.\ Jadbabaie, J.\ Lin and A.S.\ Morse (2003). 
Coordination of groups of
mobile autonomous agents using nearest neighbor rules, \emph{IEEE 
Trans.\
Automatic Control}, Vol.\ 48 (6), pp.\ 988--1001

\bibitem{JadMot2004} A. Jadbabaie, N. Motee and M. Barahona (2004).
On the stability of the Kuramoto model of coupled nonlinear 
oscillators.
In \emph{Proc.\ of the American Control Conference}

\bibitem{Naomi} N.E.\ Leonard and E.\ Fiorelli. (2001) Virtual leaders, 
artificial
potentials and coordinated control of groups. In \emph{Proc.\ of the 
40th IEEE
Conf.\ on Decision and Control}, pp.\ 2968--2973, Orlando, FL, USA

\bibitem{Lin} J.\ Lin, A.S.\ Morse and B.D.O.\ Anderson (2003). The 
multi-agent 
rendezvous problem. In \emph{Proceedings of the 42nd IEEE Conference on 
Decision and
Control}, pp.\ 1508--1513, Maui HI, USA

\bibitem{MOR03} L.\ Moreau (2003).
Time-dependent unidirectional communication in multi-agent systems,
arXiv:math.OC/0306426

\bibitem{MOR04} L.\ Moreau (2004).
A note on leaderless communication via bidirectional and unidirectional 
time-dependent communication, {\em Proc.\ of MTNS'04}, Leuven (Belgium)

\bibitem{MOR04a} L.\ Moreau (to appear).
Stability of multi-agent systems with time-dependent communication 
links,
{\em IEEE Transactions on Automatic Control}


\bibitem{Sepulchre} R.\ Sepulchre, D.\ Paley and N.\ Leonard (2003).
Collective  motion
and oscillator synchronization. 
In \emph{Proc.\ of the Block Island 
Workshop on
Cooperative Control}

\end{thebibliography}
